\documentclass{article}

\newtheorem{theorem}{Theorem}[section]
\newtheorem{corollary}{Corollary}[section]
\newtheorem{lemma}{Lemma}[section]

\newtheorem{definition}{Definition}[section]

\newtheorem{example}{Example}[section]

\begin{document}

\vspace*{.7cm}
\begin{center}
{\large{\bf  Quaternion CR-Submanifolds of a Locally Conformal
Quaternion Kaehler Manifold} }\\
\vspace{.3cm}
 {Bayram \d{S}ahin}\\
 \vspace{.3cm}
 {\it {\small {Dedicated to the memory of Dr. \"{U}mran Pekmen}}}
\end{center}
\vspace{.7cm}

\noindent {\small {\bf Abstract.{\sl  The purpose of the present
paper is to study the differential geometric properties of a
quaternion CR-submanifold in a locally
conformal quaternion Kaehler manifold.}}} \\


\section{Introduction}

The concept of the locally conformal Kaehler manifolds was
introduced by I.Vaisman in \cite{Vaisman1}. Since then many papers
appeared on these manifolds and their submanifolds ( See:
\cite{Dragomir-Ornea} and its references). However,the geometry of
the locally conformal quaternion Kaehler manifolds has been
studied in the last ten years, \cite{Dragomir-Ornea},
\cite{Ornea-Piccini},\cite{Ornea}, \cite{Pedersen1}
,\cite{Pedersen2}, \cite{Piccini1}, \cite{Piccini2} and  their
QR-Submanifolds have been studied in \cite{Sahin-Gunes}
\bigskip

A locally conformal quaternion Kaehler manifold (Shortly, l.c.q.K
manifold) is a quaternion Hermitian manifold whose metric is
conformal to a quaternion Kaehler metric in some neighborhood of
each point. The main difference between locally conformal Kaehler
manifold and l.c.q.K. manifold is that the Lee form of a compact
l.c.q.K. manifold can be choosen as parallel form without any
restrictions\cite{Dragomir-Ornea}. It is known that this property
is not guaranteed in the complex case,\cite{Tricerri}.
\bigskip

On the other hand,A.Bejancu \cite{BejancuCR} defined and studied
CR-submanifolds of a Kaehler manifold. Since then many papers
appeared on this topic,
\cite{BejancuCR},\cite{Blair-Chen},\cite{Chen3},\cite{Chen1}.
Moreover, Barros, Chen and Urbano defined quaternion CR-submanifold
of a quaternion Kaehler manifold as analogy with CR-submanifold of a
Kaehler manifold \cite{Barros-Chen-Urbano}.\\

In this paper, we study the geometry of quaternion CR-submanifolds
of a l.c.q.K. manifold.In section 2, we give some basic
definitions, formulae and result which will be used in this paper.
In section 3, we study the geometry of quaternion CR-submanifolds
of a l.c.q.K. manifold. In this section, we investigate the
geometry of leaves. In section 4, we consider totally umbilical
quaternion CR-submanifolds.  Finally, section 5 is devoted to
investigate the topology of a quaternion CR-submanifold.

\section{Locally Conformal Quaternion Kaehler Manifolds}
\setcounter{equation}{0}
\renewcommand{\theequation}{2.\arabic{equation}}

We denote a quaternion Hermitian manifold by $\left(
\bar{M},g,H\right) ,$ where $H$ is a subbundle of $End(T \bar{M})$
of rank 3 which is spanned by almost complex structures
$J_{1},J_{2},$ and $J_{3}.$ We recall that a quaternion Hermitian
metric $g$ is said to be a quaternion Kaehler metric if its
Levi-Civita connection $\bar{\nabla}$ satisfies $\bar{\nabla
}H\subset H.$

A quaternion Hermitian manifold with metric $g$ is a l.c.q.K.
manifold if over neighborhoods $\left\{ U_{i}\right\}  $ covering
$M,$ $g|_{U_{i}}=e^{f_{i}}g_{i}^{^{\prime }}$ with
$g_{i}^{^{\prime}}$ a quaternion Kaehler metric on $U_{i}.$In this
case, the Lee form $\omega$ is locally defined by
$\omega|_{U_{i}}=df_{i}$ and satisfies

\begin{equation}
d\Theta=\omega\wedge\Theta,d\omega=0 \label{bir}
\end{equation}
where $\Theta=\sum_{\alpha=1}^3\Omega_{\alpha}
\wedge\Omega_{\alpha}$ is the Kaehler $4-$form. We note that
property (\ref{bir}) is also a sufficient condition for a quaternion
Hermitian metric to be a l.c.q.K. metric \cite{Ornea-Piccini}.

Let $\bar\nabla $ be the Levi-Civita connection of $g$. We define
a Lee vector field $B$ by $\omega(X)=g(X,B)$ on $\bar{M}$. We have
another  torsionless linear connection $\bar{D}$ called the Weyl
connection. The Weyl connection $\bar{D}$ related to the
Levi-Civita connection $\bar{\nabla}$ of $g$ by the formula

\begin{equation}
\bar{D}_{X}Y=\bar{\nabla}_{X}Y-\frac{1}{2}\left\{ \omega\left(
X\right)  Y+\omega\left(  Y\right)  X-g(X,Y)B\right\} \label{iki}
\end{equation}
for any $X,Y\in\Gamma\left(  T\bar{M}\right)  ,$ where
$B=\omega^{\#}$ is the Lee vector field \cite{Dragomir-Ornea}.

Let $\bar{M}$ be a l.c.q.K. manifold and $\bar{\nabla}$ be the
connection of $\bar{M}$. Then the Weyl connection does not
preserve the compatible almost complex structures individually but
only their 3-dimensional bundle $H$, that is,
\begin{equation}
\bar{D}J_{a}=\sum Q_{ab}\otimes J_{b}  \label{uc}%
\end{equation}
for $a,b=1,2,3$,  and $Q_{ab}$ is a skew-symmetric matrix of local
forms
\cite{Pedersen1}. Thus, from (\ref{bir}) and (\ref{iki}) we have%
\begin{equation}
\begin{array}{cc}
  \bar{\nabla}_{X}J_{a}Y &= J_{a}\bar{\nabla}_{X}Y+\frac12\left\{ \theta_{o}\left(Y\right)  X \omega\left(  Y \right)  J_{a}X-\Omega left(X,Y \right)  B\\
                         & + g(X,Y)J_{a}B \} +Q_{ab}(X)J_{b}Y+Q_{ac}(X)J_{c}Y
\end{array}
\label{dort}
\end{equation}
for any $X,Y\in\Gamma\left(  T\bar{M}\right)  $, where
$\theta_{o}=\omega oJ_{a}$.

Let $R^{\bar{D}}$ and $\bar{R}$ be the curvature tensor fields of
the connections $\bar{D}$ and $\bar{\nabla}$, respectively. Then
we have
\begin{equation}
\begin{array}{cc}
  R^{\bar{D}}(X,Y)Z &=\bar{R}(X,Y)Z-\frac12\{[(\bar{\nabla}_{X}\omega)Z+\frac12\omega(X)\omega(Z)]Y \\
                          &[(\bar{\nabla}_{Y}\omega)Z+\frac12\omega(Y)\omega(Z)]X +((\bar{\nabla}_{X}\omega)Y)Z\\
                          &-((\bar{\nabla}_{Y}\omega)X)Z -g(Y,Z)(\bar{\nabla}_{X}B+\frac12\omega(X)B)\\
                          &+g(X,Z)(\bar{\nabla}_{Y}B+\frac12\omega(Y)B\}\\
                          &-\frac14\|\omega\|^2(g(Y,Z)X-g(X,Z)Y)
\end{array}
\label{bes}
\end{equation}
for any $X,Y\in\Gamma(T\bar{M})$

Next we give the following  theorem which will be useful  later.
\begin{theorem}\cite{Ornea}
Let $(\bar{M},\bar{g},H)$ be a compact quaternion Hermitian Weyl
manifold, non-quaternion Kaehler, whose foliation $F$ has compact
leaves. Then the leaves space $P=\bar{M}/{F}$ is a compact
quaternion Kaehler orbifold with positive scalar curvature, the
projection is a Riemannian, totally geodesic submersion and a
fibre bundle map with fibres as described in proposition 4.10 of
\cite{Ornea}, where $F$ is locally generated by
$B,J_{1}B=B_{1},B_{2},B_{3}$.
\end{theorem}

If $F$ is a regular foliation, then $P=\bar{M}/{F}$ is a compact
quaternion Kaehler manifold \cite{Molino}.

Let $\bar{M}$ be a l.c.q.K. manifold and $M$ be any submanifold of
$\bar{M}$. The formulae of  Gauss and Weingarten are given by

\begin{equation}
\bar{\nabla}_{X}Y=\nabla_{X}Y+h(X,Y) \label{alti}%
\end{equation}
and

\begin{equation}
\bar{\nabla}_{X}V=-A_{V}X+\nabla_{X}^{\perp}V \label{yedi}
\end{equation} for  vector fields $X,Y$ tangent to $M$ and any
vector field $V$ normal to $M$, where $\nabla$ is the induced
Riemann connection in $M,$ $h$ is the second fundamental form,
$A_{V}$is the fundamental tensor field of Weingarten with respect to
the normal section $V$ and $\nabla^{\perp}$ is the
normal connection. Moreover, we have the relation%

\begin{equation}
g(h(X,Y),V)=g(A_{V}X,Y). \label{sekiz}
\end{equation}

The equations of Gauss and  Codazzi are given respectively  by
\cite{Chen}

\begin{equation}
\begin{array}{cc}
  R(X,Y;Z,W) &= \bar{R}(X,Y;Z,W)+g(h(X,W),h(Y,Z)) \\
             & -g(h(X,Z),h(Y,W))
\end{array}
\label{Gauss}
\end{equation}
and
\begin{equation}
(\bar{R}(X,Y)Z)^\perp=(\bar{\nabla}_{X}h)(Y,Z)-(\bar{\nabla}_{Y}h)(X,Z)
\label{Codazzi}
\end{equation}
for $X,Y,Z\in\Gamma(TM)$, where $\bar{R},R$ is the curvature
tensor corresponding to the connection $\bar{\nabla},\nabla$
respectively and $ ^\perp$ in (\ref{Codazzi}) denotes the normal
component. The covariant derivative $(\bar{\nabla}_{X}h)(Y,Z)$ is
given by
$$(\bar{\nabla}_{X}h)(Y,Z))=\nabla^{\perp}_{X}h(Y,Z)-h(\nabla_{X}Y,Z)-h(Y,\nabla_{X}Z).$$

\section{Quaternion CR-Submanifolds of a l.c.q.K. manifold}
\setcounter{equation}{0}
\renewcommand{\theequation}{3.\arabic{equation}}

First, we give definition of a quaternion CR-submanifold of a
l.c.q.K. manifold as the definition of quaternion CR-submanifolds
of a quaternion Kaehler manifold.
\begin{definition} A submanifold $M$ of a l.c.q.K. manifold
$\bar{M}$ is called a quaternion CR-submanifold if there exists
two orthogonal complementary distributions $D$ and $D^\perp$  such
that $D$ is invariant under $J_{a}$,i.e,$J_{a}D\subseteq D$,
$a=1,2,3.$ and $D^\perp$ is totally real,i.e. $J_{a}D^\perp
\subseteq TM^\perp$, $ a=1,2,3.$
\end{definition}

A submanifold $M$ of a l.c.q.K. manifold $\bar{M}$ is called a
quaternion submanifold(resp. totally real submanifold) if $dim
D^\perp=0$ (resp. dim$D$=0). A quaternion CR-submanifold is called
proper quaternion CR-submanifold if it is neither quaternion nor
totally real.

By the definition a quaternion CR-submanifold, we have
\begin{equation}
TM=D\oplus D^\perp \label{dokuz}
\end{equation}
and
\begin{equation}
TM^\perp=J_{a}D^\perp\oplus\mu \label{on}
\end{equation}
where $\mu$ is orthogonal complement of $J_{a}D^\perp$ in the
normal bundle is invariant subbundle of  $\Gamma(TM^\perp)$ under
$J_{a}$.

Now, let $M$ be a quaternion CR-submanifold of a l.c.q.K. manifold
$\bar{M}$. For each vector field  $X$ tangent to $M$ we put
\begin{equation}
J_{a}X=\phi_{a}X+\varpi_{a}X \label{onbir}
\end{equation}
where $\phi_{a}X\in\Gamma(D)$ and
$\varpi_{a}X\in\Gamma(J_{a}D^\perp)$. Also, for  each vector field
$V$ normal to $M$ we put
\begin{equation}
J_{a}V= f_{a}V+t_{a}V \label{oniki}
\end{equation}
where $f_{a}V\in\Gamma(D^\perp)$ and $t_{a}V\in\Gamma(\mu)$.

Now, we will  give an example for quaternion CR-submanifolds of a
l.c.q.K. manifold.
\begin{example}
Let $\bar{M}$ be a l.c.q.K. manifold. Assume that the foliation
$F$ is regular. Then $P=\bar{M}/{F}$ is a compact quaternion
Kaehler manifold (cf:Theorem.2.1). We denote almost complex
structures of $\bar{M}$ and $P$ by $J_{a}$ and $J'_{a}$,
respectively. Now we consider the following commutative diagram:
\[
   \begin{array}{ccc}
\bar{M}                                 &
\stackrel{\pi}{\longrightarrow}              &
 P=\bar{M}/{F}                                     \\
\Big\uparrow \vcenter{%
   \rlap{$ \scriptstyle{i}$}}                   & &
\Big\uparrow \vcenter{%
\rlap{$ \scriptstyle{j}$}}                          \\
 N                                                &
\stackrel{\bar{\pi}} {\longrightarrow}
                                           & \bar{N}
\end{array}
\]
where $N$ and $\bar{N}$ are submanifolds of $\bar{M}$ and $P$,
respectively.We denote  the horizontal lift by $^*$. Then we have
\begin{equation}
(J'_{a}X)^*=J_{a}X^*. \label{sub1}
\end{equation}
We note that the projection $\pi$ is a totally geodesic Riemannian
submersion and a fibre bundle map. Hence $\bar{\pi}$ is also a
Riemannian submersion. We denote the vertical distribution of the
Riemannian submersion $\pi$ by $\upsilon$, i.e.
$ker\pi_{*}=\upsilon$. Let $\bar{H}$ be  the horizontal
distribution of $\pi$. Then we have
$T\bar{M}=\bar{H}\oplus\upsilon$. We denote the horizontal
distribution of $\bar{\pi}$ by $H_{0}$.We will investigate the
relation between  normal spaces of $N$ and $\bar{N}$. We denote
the  Riemannian metrics of $\bar{M}$ and $P$ by $g$ and $g'$,
respectively. Let $V^*$ be the horizontal lift of
$V\in\Gamma(T\bar{N}^\perp)$. Then we get
\begin{eqnarray*}
g(V^*,X)&=&g((\pi_{*})^*V,X)\\
&=&g'(\pi_{*}X,V) \\
&=&0,
\end{eqnarray*}
for any $X \in H_{0}$. Thus, $(T\bar{N}^\perp)^*$ is orthogonal to
$H_{0}$. Note that the normal space is always horizontal. Hence
$(T\bar{N}^\perp)^*$ is orthogonal to ${\upsilon}$. Consequently,
we have $(T\bar{N}^\perp)^* \subseteq TN^\perp$. Since $\pi$ is a
Riemannian submersion we get

\begin{equation}
(T\bar{N}^\perp)^*=TN^\perp. \label{sub2}
\end{equation}

 Now, let $\phi_{a}$ and $\varpi_{a}$ be  the operators on $\bar{N}$
 appearing in (\ref{onbir}). We denote  the  operators  in $N$ corresponding  to $\phi_{a}$ and $\varpi_{a}$ by $\phi'_{a}$ and $\varpi'_{a}$, respectively. From (\ref{sub1}) and (\ref{sub2}) we obtain
 \begin{equation}
 (\phi_{a}X)^*=\phi'_{a}X^* \label{onuc}
 \end{equation}
 and
 \begin{equation}
 (\varpi_{a}X)^*=\varpi'_{a}X^*. \label{ondort}
 \end{equation}
 So, from  (\ref{onuc}) and (\ref{ondort}) we see that  $N$ is
 a quaternion CR-submanifold of $\bar{M}$ if and only if $\bar{N}$
 is a quaternion CR-submanifold of $P$.
\end{example}
In the  rest of this section, we will investigate  the geometry of
leaves on quaternion CR-submanifolds.

In \cite{Barros-Chen-Urbano}, Barros, Chen and Urbano showed that
the anti-invariant distribution $D^\perp$ of a quaternion Kaehler
manifold is integrable. In the next theorem, we will see that is
still true for a quaternion CR-submanifold of a l.c.q.K. manifold.
\begin{theorem} Let $M$ be a proper quaternion CR-submanifold of a
l.c.q.K. manifold. Then the anti-invariant distribution $D^\perp$
of $M$ is integrable.
\end{theorem}
\noindent{\bf Proof.~} From (\ref{dort}),(\ref{alti}) and
(\ref{yedi}) we obtain
\begin{equation}
-g(A_{J_{a}W}T,X)=g(\nabla_{T}W,X)+g(T,W)g(J_{a}B,X) \label{ondokuz}
\end{equation}
for any $X\in\Gamma(D)$ and $T,W\in\Gamma(D^\perp)$.Interchanging
$T$ and $W$ in (\ref{ondokuz}) and subtracting we get
\begin{equation}
g(A_{J_{a}T}W-A_{J_{a}W}T,X)=g(J_{a}[T,W],X). \label{yirmi}
\end{equation}
On the other hand, since $\bar{\nabla}$ is a metric connection and
$A$ is self-adjoint we obtain
\begin{equation}
g(A_{J_{a}T}W,X)=-g(W,\bar{\nabla}_{X}J_{a}T).
\end{equation}
In this equation, using (\ref{dort}) and (\ref{alti}) we have
\begin{equation}
g(A_{J_{a}T}W,X)=g(A_{J_{a}W}T,X). \label{yirmibir}
\end{equation}
Thus, from (\ref{yirmi}) and (\ref{yirmibir}), we obtain
$$g([T,W],J_{a}X)=0,$$ which proves our assertion.
\\

\begin{lemma}Let $M$  be a quaternion CR-submanifold of a l.c.q.K.
manifold. Then quaternion distribution $D$ is minimal if and only
if the Lee vector field is orthogonal to the anti-invariant
distribution $D^\perp$.
\end{lemma}
\noindent{\bf Proof.~} Since $\bar{\nabla}$ is a metric connection,
from (\ref{dort}) and (\ref{yedi}) we obtain
\begin{equation}
g(\nabla_{X}X,Z)=g(A_{J_{a}Z}X,J_{a}X)-\frac12\|X\|^2\omega(Z)
\label{yirmiiki}
\end{equation}
for any $X\in\Gamma(D)$ and $Z\in\Gamma(D)$. In a similar way we
have
\begin{equation}
g(\nabla_{J_{a}X}{J_{a}X},Z)=-g(A_{J_{a}Z}X,J_{a}X)-\frac12\|X\|^2\omega(Z)
\label{yirmiuc}
\end{equation}
Thus from (\ref{yirmiiki}) and (\ref{yirmiuc}) we have
$$
g(\nabla_{X}X,Z)+g(\nabla_{J_{a}X}{J_{a}X},Z)=0\Leftrightarrow\omega(Z)=0
$$
\\

Now, we will discuss the integrability of the quaternion
distribution. First we give a lemma.
\begin{lemma}
Let $M$ be a quaternion CR-submanifold of a l.c.q.K. manifold
$\bar{M}$. Then we have
\begin{equation}
h(X,J_{a}Y)=\varpi_{a}\nabla_{X}Y+t_{a}h(X,Y)+\frac12\{g(X,Y)(J_{a}B)^\perp-\Omega(X,Y)B^\perp
\} \label{yirmidort}
\end{equation}
and
\begin{equation}
h(X,J_{a}Y)-h(Y,J_{a}X)=\varpi([X,Y])+\Omega(X,Y)B^\perp
\label{yirmibes}
\end{equation}
for any $X,Y\in\Gamma(D)$, where $B^\perp=norB$.
\end{lemma}
\noindent{\bf Proof.~} Using (\ref{dort}), (\ref{alti}),
(\ref{yedi}) and a comparison between normal components we derive
the (\ref{yirmidort}). Then (\ref{yirmibes}) is direct consequence
of (\ref{yirmidort}).
\\
\begin{definition}
Let $M$ be a quaternion CR-submanifold of a l.c.q.K. manifold
$\bar{M}$. Then $M$ is called $D$-geodesic if $h(X,Y)=0$ for
$X,Y\in\Gamma(D)$.
\end{definition}
From Lemma 3.2 and Definition.3.2 we have:
\begin{theorem}
Let $M$ be a quaternion CR-submanifold of a l.cq.K. manifold
$\bar{M}$. Assume that the Lee vector field is tangent to $M$.
Then the following assertions are equivalent:
\begin{enumerate}
  \item for any $X,Y\in\Gamma(D)$
 $$h(X,J_{a}Y)=h(Y,J_{a}X).$$
  \item $M$ is $D$-geodesic.
  \item $D$ is integrable.
\end{enumerate}
\end{theorem}

The proof is similar to that of Theorem 2.1 in \cite{Bejancu}
(Also, Theorem 3.2 in \cite{Sahin-Gunes}). So we omit it here.

From Lemma 3.2 we have the following corollary.
\begin{corollary}
Let $M$ be a quaternion CR-submanifold of a l.cq.K. manifold
$\bar{M}$.If $M$ is $D$-geodesic and the Lee vector field
orthogonal to $D^\perp$, then each leaf of $D$ is totally
geodesic.
\end{corollary}
\begin{corollary}
Let $M$ be a quaternion CR-submanifold of a l.cq.K. manifold
$\bar{M}$. If $h(X,Z)\subset\mu$ for $X\in\Gamma(D)$
,$Z\in\Gamma(D^\perp)$ and the Lee vector field is orthogonal to
$D$, then each leaf of $D^\perp$ is totally geodesic in $M$.
\end{corollary}
\noindent{\bf Proof.~}From (\ref{dort}), (\ref{alti}) and
(\ref{yedi}) we obtain
$$
-g(A_{J_{a}W}Z,X)=-g(\nabla_{Z}W,J_{a}X)+\frac12g(Z,W)g(J_{a}B,X)
$$
for any $Z,W\in\Gamma(D^\perp)$ and $X\in\Gamma(D)$. Using
(\ref{sekiz}) we arrive at
$$
g(h(Z,X),J_{a}W)=g(\nabla_{Z}W,J_{a}X)+\frac12g(Z,W)g(B,J_{a}X)
$$
which proves our assertion.
\\
\section{Umbilical Quaternion CR-Submanifolds of l.c.q.K.
Manifolds} \setcounter{equation}{0}
\renewcommand{\theequation}{4.\arabic{equation}}

Let $\bar{M}$ be a compact l.c.q.K. manifold. Then we can choose
the fixed metric $g$ such that

i) The fixed metric $g$ makes $\omega$ parallel, i.e.
\begin{equation}
\bar{\nabla}\omega=0 \label{yirmialti}
\end{equation}
ii)
\begin{equation}
\|\omega\|^2=1 \label{yirmiyedi}
\end{equation}
(See:\cite{Ornea-Piccini}). From now on we will denote a compact
l.c.q.K. manifold by $\bar{M}$ in this section.

Let $M$ be a quaternion CR-submanifold of l.c.q.K. manifold
$\bar{M}$. Then $M$ is called totally umbilical if we have
\begin{equation}
h(X,Y)=g(X,Y)H \label{yirmisekiz}
\end{equation}
for any $X,Y$ tangent to $M$, where $H$ is the mean curvature
vector field defined by $ H=\frac1m Trace(h)$. We say that $M$ is
totally geodesic if $h=0$ identically on $M$.
\begin{theorem}
Let $M$ be a quaternion CR-submanifold of a l.c.q.K. manifold If
the Lee vector field $B$ is tangent to $M$, then we have
\begin{equation}
K_{\bar{M}}(X,Y)\leq\frac{1}{4} \label{otuzdort}
\end{equation}
for any orthonormal vector fields $X\in\Gamma(D)$ and
$Y\in\Gamma(D^\perp)$. The equality (\ref{otuzdort}) holds if and
only if $\bar{M}$ is a quaternion Kaehler manifold
\end{theorem}
\noindent{\bf Proof.~}
 Let $R^{\bar{D}}$ be the
curvature tensor field of the Weyl connection $\bar{D}$. Then we
have
\begin{equation}
R^{\bar{D}}(X,Y)J_{a}Z-J_{a}R^{\bar{D}}(X,Y)Z=\alpha(X,Y)J_{b}Z-\beta(X,Y)J_{c}Z
\label{yirmidokuz}
\end{equation}
where
$$ \alpha=dQ_{ab}+Q_{cb}\wedge Q_{ac}$$
and
$$ \beta=dQ_{ac}+Q_{bc}\wedge Q_{ab}.$$
Taking account (\ref{bes}), (\ref{yirmialti}), (\ref{yirmiyedi})
and (\ref{yirmidokuz}) we obtain
\begin{equation}
 \begin{array}{cc}
   0 &= -\bar{R}(X,Y,X,Y)+\bar{R}(X,Y,J_{a}X,J_{a}Y) \\
     & \frac{1}{4}\omega(X)\omega(X)+\frac{1}{4}\omega(Y)\omega(Y)-\frac{1}{4}
 \end{array}
\label{otuz}
\end{equation}
for any orthonormal vector fields $X\in\Gamma(D)$ and
$Y\in\Gamma(D^\perp)$. On the other hand, from (\ref{yirmisekiz})
and Codazzi equation we have
\begin{equation}
\bar{R}(X,Y,Z,V)=g(Y,Z)g(\nabla^\perp_{X}H,V)-g(X,Z)g(\nabla^\perp_{Y}H,V)
\label{otuzbir}
\end{equation}
for any vector fields $X,Y,Z$ tangent to $M$ and $V$ normal to
$M$. Thus using (\ref{otuzbir}) we get
\begin{equation}
\bar{R}(X,Y,J_{a}X,J_{a}Y)=0 \label{otuziki}
\end{equation}
for any $X\in\Gamma(D)$ and $Y\in\Gamma(D^\perp)$. Using
(\ref{otuz}) and (\ref{otuziki}) we arrive at
\begin{equation}
\bar{R}(X,Y,Y,X)=-\frac{1}{4}\omega(X)^2-\frac{1}{4}\omega(Y)^2+\frac{1}{4}.
\label{otuzuc}
\end{equation}

If the Lee vector field is tangent to $M$  from (\ref{otuzuc}) we
have (\ref{otuzdort}). In view of (\ref{otuzuc}) the equality case
of (\ref{otuzdort}) is valid if and only if
$$\omega(X)^2+\omega(Y)^2=0$$
for any orthonormal vector fields $X\in\Gamma(D)$ and
$Y\in\Gamma(D^\perp)$.Thus we obtain
\begin{equation}
\omega(X)=0 \label{otuzbes}
\end{equation}
and
\begin{equation}
\omega(Y)=0. \label{otuzalti}
\end{equation}
From (\ref{otuzbes}) and (\ref{otuzalti}) we obtain $B\in\Gamma(D)$
and $B\in\Gamma(D^\perp)$, respectively. Since $D \cap
D^\perp=\{0\}$ we have $B=0$.\\

Now we give an another   theorem  for totally umbilical proper
quaternion CR-submanifolds of a l.c.q.K. manifold. We start with
the following preparatory result.
\begin{lemma}
Let $M$ be a totally umbilical quaternion CR-submanifold of a
l.c.q.K. manifold. Assume that the Lee vector field is tangent to
$M$. Then we have
\begin{equation}
H\in\Gamma(J_{a}D^\perp) \label{otuzdokuz}
\end{equation}
\end{lemma}
\noindent{\bf Proof.~} Since $B$ is tangent to $M$, from Lemma 3.2
we obtain
$$g(h(X,J_{a}Y),N)=-g(h(X,Y),J_{a}N)$$
for any $N\in\Gamma(\mu)$ and $X,Y\in\Gamma(D)$. Since $M$ is
totally umbilical we get
$$g(X,J_{a}Y)g(H,N)=-g(X,Y)g(H,J_{a}N).$$
Thus for $X=J_{a}Y$ we have $g(H,N)=0$, hence we obtain
$H\in\Gamma(J_{a}D^\perp)$
\\
\begin{theorem}
Let $M$ be a  totally umbilical proper quaternion CR-submanifold
of a l.c.q.K. manifold. Assume that the Lee vector field is
tangent to $M$. Then
\begin{enumerate}
  \item $M$ is totally geodesic or
  \item the totally real distribution is one dimensional.
\end{enumerate}
\end{theorem}
\noindent{\bf Proof.~} We take $Z,W\in\Gamma(D)$ and using totally
umbilicalness of $M$ together (\ref{dort}), (\ref{alti}) and
(\ref{yedi}) we have
$$-A_{J_{a}W}Z=\phi_{a}\nabla_{Z}W+g(Z,W)J_{a}H+\frac{1}{2}\theta(W)Z+\frac{1}{2}g(Z,W)(J_{a}B)^T.$$
Taking inner product with $Z$ in $D^\perp$ it follows that
$$
 -g(A_{J_{a}W}Z,Z)=-g(Z,W)g(H,J_{a}Z)+\frac{1}{2}\theta(W)g(Z,Z)-\frac{1}{2}g(Z,W)g(B,J_{a}Z)
$$
Since $B$ is tangent to $M$, from (\ref{sekiz}) and
(\ref{yirmisekiz}) we have
\begin{equation}
g(Z,Z)g(H,J_{a}W)=\frac{1}{2}g(Z,W)g(H,J_{a}Z). \label{kirk}
\end{equation}
Interchanging $Z$ and $W$ in (\ref{kirk}) we obtain
\begin{equation}
g(W,W)g(H,J_{a}Z)=\frac{1}{2}g(Z,W)g(H,J_{a}W). \label{kirkbir}
\end{equation}
From (\ref{kirk}) and (\ref{kirkbir}), one can immediately get
\begin{equation}
g(H,J_{a}W)=\frac{g(Z,W)^2}{\parallel Z \parallel^2 \parallel W
\parallel^2}g(H,J_{a}W). \label{kirkiki}
\end{equation}
From Lemma 4.1., the possible solutions of (\ref{kirkiki}) are \\
a) $H=0$ b)$ Z//W$\\
Suppose condition a) holds,i.e., $H=0$ which implies that $M$ is
totally geodesic. If b) is satisfies in (\ref{kirkiki}) then
dim$D^\perp=1$ which implies that the totally real distribution is
one dimensional.
\\

\section{Cohomology of Quaternion CR-Submanifolds}
\setcounter{equation}{0}
\renewcommand{\theequation}{5.\arabic{equation}}

Assume that $M$ be a quaternion CR-submanifold of $4n$ dimensional
l.c.q.K. manifold $\bar{M}$. Let $p=dim_{Q}D, q=dimD^\perp$ Then we
choose an orthonormal frame $ \{e_{1},...,e_{p}$,
$e_{p+1}=J_{1}e_{1},...,e_{4p}$, $E_{1},...,E_{q}$,
$J_{1}E_{1},...,J_{1}E_{q}$, $J_{2}E_{1},...,J_{2}E_{q}$,
$J_{3}E_{1},...,J_{3}E_{q}$, $V_{1},...,V_{r}$,
$V_{r+1}=J_{1}V_{1},...,V_{4r}\}$ in $\bar{M}$ such that restricted
to $M$,$\{e_{1},...,e_{p},e_{p+1}=J_{1}e_{1},...,e_{4p}\}$ are in
$D$ and $\{E_{1},...,E_{q}\}$ are in $D^\perp$. We denote by
$\{w^1,...,w^4p\}$ the 1-forms on $M$ satisfying
\begin{equation}
\begin{array}{ccc}
  w^i(Z)=0 &w^i(e_{j})=\delta_{ij} &i,j=1,...,4p
\end{array}
\label{ellibir}
\end{equation}
for any $Z\in\Gamma(D^\perp)$, where
$e_{p+j}=J_{1}e_{j},e_{2p+j}=jJ_{2}e_{j},e_{3p+j}=J_{3}e_{j}$Then
\begin{equation}
w=w^1 \wedge...\wedge w^4p   \label{elliiki}
\end{equation}
defines a $4p$-form on $M$. From (\ref{elliiki}) we obtain
\begin{equation}
dw=\sum_{i=1}^{4p} (-1)^i w^1 \wedge...\wedge...\wedge w^{4p}
\label{elliuc}
\end{equation}
Thus from (\ref{ellibir}) and (\ref{elliuc}) we obtain that $dw=0$
if and only if
\begin{equation}
dw(Z_{1},Z_{2},X_{1},...,X_{4p-1})=0 \label{ellidort}
\end{equation}
and
\begin{equation}
dw(Z_{1},X_{1},...,X_{4p})=0 \label{ellibes}
\end{equation}
for any $Z_{1},Z_{2}\in\Gamma(D^\perp)$ and
$X_{1},...,X_{4p}\in\Gamma(D)$. We see that (\ref{ellidort}) holds
if and only if $D^\perp$ is integrable and (\ref{ellibes}) holds
if and only if $D$ is minimal.Thus from Theorem.3.1 and Lemma.3.1,
we have the following theorem.
\begin{theorem} Let $M$ be a closed quaternion CR-submanifold of a
l.c.q.K. manifold. If the Lee vector field is orthogonal to the
anti-invariant distribution $D^\perp$ then the $4p$-form $w$
defines a canonical de Rham cohomology class [$w$] in
$H^{4p}(M,R)$
\end{theorem}

\noindent Inonu University\\
\noindent Faculty of Science and Art\\
\noindent Department of Mathematics\\
\noindent 44280-Malatya/TURKEY.\\
\noindent E-mail:bsahin@inonu.edu.tr


\begin{thebibliography}{99}



\bibitem{Barros-Chen-Urbano}M.Barros,B.Y.Chen and F.Urbano,
Quaternion CR-Submanifolds of Quaternion Kaehler Manifolds, Kodai
Math. 4 (1981) 399-418.

\bibitem {BejancuCR}A. Bejancu, CR-Submanifolds of a Kaehler
Manifold, I-II, Proc. Amer. Math. Soc., Vol:69 (1978),135-142,
Trans.Amer.Math.Soc.,250 (1979),333-345.

\bibitem {Bejancu}A. Bejancu, QR-submanifolds of Quaternion Kaehlerian
Manifolds, Chinese J. Math, Vol:14, No:2, (1986).

\bibitem{Blair-Chen} D.E.Blair and B.Y.Chen, On CR-Submanifolds
of Hermitian Manifolds, Israel J.Math.34 (1979), 353-363.

\bibitem{Chen} B.-Y. Chen, Geometry of Submanifolds, Marcel-Dekker
Inc., (1973).

\bibitem{Chen3} B.Y.Chen, CR-Submanifolds of a Kaehler Manifold,
I-II, J.Diff.Geometry,16 (1981), 305-322, 493-509.

\bibitem{Chen1} B.Y.Chen, Cohomology of CR-Submanifolds, An.
Faculte des Sciences Toulouse 3 (1981), 167-172.

\bibitem {Dragomir-Ornea}S. Dragomir, L. Ornea, Locally Conformal K\"{a}hler Geometry, Birkh\"{a}user,
(1998).

\bibitem {Matsumoto}K. Matsumoto, On CR-Submanifolds of Locally Conformal
Kaehler Manifolds, J.Korean Math. Soc., 21, No:1, (1984), 49-61.

\bibitem{Molino} P.Molino, Riemannian Foliations, Birkhauser
(1988).


\bibitem {Ornea-Piccini}L. Ornea and P. Piccini, Locally Conformal K\"{a}hler Structures in Quaternionic Geometry,
Trans. Amer. Math. Soc., Vol:349, No:2,  (1997), 641-645.

\bibitem {Ornea}L. Ornea, Weyl Structures on Quaternionic Manifolds. A State of the
Art, arXiv:math. DG/0105041, v1, (2002).

\bibitem {Pedersen1}H. Pedersen, Y.S. Poon, A. Swann, The Einstein-Weyl
Equations in Complex and Quaternionic Geometry, Diff. Geo. and Its
Appl.,  (1993), 309-322.

\bibitem {Pedersen2}H. Pedersen, A. Swann, Riemannian Submersions, Four
Manifolds and Einstein Weyl Geometry, Proc. London Math. Soc.,
(3)66, (1993), 338-351.

\bibitem{Piccini1}P.Piccini, Manifolds with local Quaternion
Kaehler Structures,Rend.Mat.17 (1997),679-696.

\bibitem{Piccini2}P.Piccini, The Geometry of Positive Locally
Quaternion Kaehler Manifold, Ann.Global Anal.Geom. 16 (1998).

\bibitem{Sahin-Gunes} B.Sahin and R.Gunes, QR-Submanifolds of a
Locally Conformal Quaternion Kaehler Manifolds, Publ.Math.Debrecen
2003, Vol:63.

\bibitem{Tricerri}F.Tricerri, Some Examples of locally Conformal
Kaehler Manifolds, Rend.Sem.Mat.Univ. Politecn.Torino, 40 (1981),
81-92

\bibitem {Vaisman1}I. Vaisman, On Locally Conformal Almost K\"{a}hler Manifolds
Israel J. Math., 24, (1976), 338-351,.

\bibitem {Vaisman2}I. Vaisman, A Theorem on Compact Locally Conformal K\"{a}hler
Manifolds, Proc. Amer. Math. Soc., Vol:75, No:2, (1979,279-283.

\bibitem {Vaisman3}I. Vaisman, On Locally and Globally Conformal K\"{a}hler
Manifolds, Trans. Amer. Math. Soc., Vol:262, No:2, (1980),439-447.

\bibitem {Vaisman4}I. Vaisman, Some Curvature Properties of Locally Conformal
K\"{a}hler Manifolds, Trans. Amer. Math. Soc., Vol:259, No:2,
(1980), 439-447.

\bibitem {Vaisman5}I. Vaisman, A Geometric Condition For An l.c.K. Manifold To
Be K\"{a}hler, Geometriae Dedicata, 10, (1981, 129-134.


\bibitem{Yano} K.Yano-M.Kon, Structures on Manifolds, Ser.Pure
Math.World Scientific, (1984)
\end{thebibliography}
\end{document}